\documentclass{amsart}

\newtheorem{theorem}{Theorem}[section]
\newtheorem{lemma}[theorem]{Lemma}

\newtheorem{cor}[theorem]{Corollary}
\theoremstyle{definition}
\newtheorem{definition}[theorem]{Definition}

\newtheorem{conjecture}[theorem]{Conjecture}

\newcommand{\BB}{\mathbb{B}}
\newcommand{\ZZ}{\mathbb{Z}}

\newcommand{\QQ}{\mathbb{Q}}

\newcommand{\RR}{\mathbb{R}}

\def\vol{\mbox{\rm{Vol}}}
\def\Area{\mbox{\rm{Area}}}

\def\T{\Theta}

\def\D{\Delta}

\def\disc{\mathcal{D}}
\def\deg{\mbox{\rm{deg}}}
\def\cO{{\mathcal O}}

\def\e{\epsilon}
\def\g{\gamma}

\def\G{\Gamma}
\def\th{\theta}
\def\l{\lambda}

\def\s{\sigma}

\def\i{\iota}

\def\S{\Sigma}

\def\ie{{\it i.e.}\ }

\def\tr{\mbox{\rm{tr}}}

\def\PSL{\mbox{\rm{PSL}}}

\def\Isom{\mbox{Isom}}
\def\Mob{\mbox{M\"ob}}

\def\GL{\mbox{\rm{GL}}}

\def\O{\mbox{\rm{O}}}
\def\OO{\mathcal{O}}

\def\Vol{\mbox{\rm{Vol}}}
\def\Comm{\mbox{\rm{Comm}}}

\def\mod{\mbox{mod}}

\def\HH{\mathbb{H}}
\def\SS{\mathbb{S}}

\def\EE{\mathbb{E}}

\def\interior{\rm int}

\numberwithin{equation}{section}

\begin{document}

% \title[short text for running head]{full title}
\title{Finiteness of arithmetic hyperbolic reflection groups}

%    author one information
% \author[short version for running head]{name for top of paper}
\author{Ian Agol}
\address{MSCS UIC 322 SEO, m/c 249\\
        851 S. Morgan St.\\
        Chicago, IL 60607-7045\\
        USA}
\email{agol@math.uic.edu}
\thanks{Agol partially supported by NSF grant DMS-0504975 and the Guggenheim Foundation}

%    author two information
\author{Mikhail Belolipetsky}
\address{Department of Mathematical Sciences\\
Durham University\\
South Rd\\
Durham, DH1 3LE\\
UK}
\email{mikhail.belolipetsky@durham.ac.uk}
\thanks{Belolipetsky thanks MPIM in Bonn for hospitality and
        support during the summer of 2006}

\author{Peter Storm}
\address{Stanford University\\
Mathematics, Bldg. 380\\
450 Serra Mall\\
Stanford, CA 94305-2125\\
USA}
\email{storm@math.stanford.edu}
\thanks{Storm partially supported
        by NSF grant DMS-0603711}

\author{Kevin Whyte}
\address{MSCS UIC 322 SEO, m/c 249\\
        851 S. Morgan St.\\
        Chicago, IL 60607-7045\\
        USA}
\email{kwhyte@math.uic.edu}
\thanks{Whyte partially funded by 
NSF Career grant DMS-0349290 and
the Sloan Foundation}

%    \subjclass is required.

\date{December 4, 2006}

%\dedicatory{}
\large
%    Abstract is required.
\begin{abstract}
We prove that there are only finitely many conjugacy classes of arithmetic maximal hyperbolic
reflection groups.
\end{abstract}

\maketitle

%    Text of article.

\section{Introduction}
A hyperbolic reflection group $\G$ is a discrete group generated by reflections in the
faces of a hyperbolic polyhedron $P\subset \HH^{n}$. We may assume that the
dihedral angles of $P$ are of the form $\pi/n$, $n\geq 2$, in which case $P$ forms
a fundamental domain for the action of $\G$ on $\HH^{n}$. If $P$ has finite volume,
then $\HH^{n}/\G=\mathcal{O}$ is a hyperbolic orbifold of finite volume, which is
obtained by ``mirroring'' the faces of $P$.
A reflection group $\G$ is maximal if there is no reflection group $\G'$ such that
$\G < \G'$. We defer the definition of arithmetic
groups until later, but a theorem of Margulis implies that $\G$ is arithmetic if and only if $[\Comm(\G):\G]=\infty$,
where $\Comm(\G)=\{g\in \Isom(\HH^{n}) | [\G:g^{-1}\G g \cap \G]<\infty\}$ \cite{Margulis75}.
The main theorem of this paper is that up to conjugacy in $\Isom(\HH^{n})$ there are only finitely many
arithmetic groups which are maximal hyperbolic reflection groups. If one does not assume that
the groups are maximal, then there is no corresponding finiteness result. For example,
using a method of \cite{S} one can see that there
are all right angle reflection groups which then contain infinitely many finite index reflection subgroups, by
a sequence of doublings along faces.
These examples are known up to dimension $8$, but recently examples
of infinite families of arithmetic reflection groups up to dimension 19
have been given by Allcock \cite{Allcock06}.

Finiteness results of this sort were proven before by Nikulin \cite{Nikulin80, Nikulin81, Nikulin87}
for dimension $n\ge 10$. Subsequently Vinberg
showed that there are no uniform (\ie cocompact) reflection groups in dimension $n\ge 30$ \cite{Vinberg84}.
Our results are independent of their work, in the sense that
we obtain a finiteness result in each dimension and then appeal to the result of Prokhorov that reflection groups
do not exist in high enough dimensions \cite{Prokhorov86}. 
The argument of this paper generalizes an argument of Long-MacLachlan-Reid \cite{LMR06}, which implies that
there are only finitely many arithmetic minimal (or congruence) hyperbolic 2-orbifolds with bounded
genus. Their argument is in fact a generalization of an argument of Zograf \cite{Zograf91}, who
reproved that there are only finitely many congruence groups $\G$ commensurable with
$\PSL(2,\ZZ)$ such that $\HH^{2}/\G$ has genus $0$
(This was proven originally by Dennin \cite{Dennin71,Dennin72}, and was known as
Rademacher's conjecture). These arguments were generalized in \cite{Agol05} to show that
there are only finitely many  maximal arithmetic reflection groups in dimension 3.
Thus, our result is only new for cocompact maximal arithmetic groups in dimensions $4\leq n \leq 9$. 
One key ingredient of these arguments is a theorem of Burger and Sarnak \cite{BurgerSarnak91},
which implies that there is a lower bound on the first eigenvalue of a congruence arithmetic
hyperbolic group defined by a quadratic form.
The other key ingredient is an inequality of Li-Yau \cite{LiYau82},
% El Soufi-Ilias \cite{ElSoufiIlias86},
generalized to orbifolds, which allows us to give an upper bound on the first eigenvalue of a maximal
quotient of an arithmetic reflection group. One of the new ingredients of the present work is
the use of some results on the distribution of covolumes of maximal arithmetic subgroups obtained
in \cite{Bel}.

{\bf Remark: } As this manuscript was being prepared, a short preprint appeared by Nikulin proving 
the finiteness of arithmetic maximal reflection groups in all dimensions \cite{Nikulin06}. Nikulin's proof makes 
use of Nikulin's previous results \cite{Nikulin80, Nikulin81}, as well as the finiteness results 
for maximal arithmetic reflection groups in dimensions two and three \cite{LMR06, Agol05}. The
methods in this paper are completely independent of Nikulin's methods. A final classification of maximal
arithmetic reflection groups might combine the techniques of both approaches, so it seemed 
interesting to publish our result, even with the appearance of Nikulin's preprint. 

\section{Conformal volume of orbifolds} \label{conformal volume}
Conformal volume was first defined by Li-Yau \cite{LiYau82}, partially motivated by generalizing results
on surfaces due to Yang-Yau \cite{YangYau80}, Hersch \cite{Hersch70}, and Szeg\H{o} \cite{Szego54}. We generalize this notion to orbifolds.
Let $(\mathcal{O},g)$ be a complete Riemannian orbifold, possibly
with boundary. Let
$|\mathcal{O}|$ denote the underlying topological space. Denote
the volume form by $dv_{g}$, and the volume by $\Vol(\mathcal{O},g)$. Let $\Mob(\SS^{n})$ denote
the conformal transformations of $\SS^{n}$. It is well-known that $\Mob(\SS^{n})=\Isom(\HH^{n+1})$.
The topological space $|\mathcal{O}|$ has a dense open subset which is a Riemannian manifold. We will
say that a map
$\varphi: |\mathcal{O}_{1}|  \to |\mathcal{O}_{2}|$ is {\it PC} if it is a continous map which is piecewise a conformal immersion.
Clearly, if $\varphi:|\mathcal{O}|\to \SS^{n}$ is PC, and $\mu\in \Mob(\SS^{n})$, then $\mu\circ \varphi$ is also PC.
Let $(\SS^{n},can)$ be the canonical round metric on the $n$-sphere.

\begin{definition}
For a piecewise smooth map $\varphi: |\mathcal{O}| \to (\SS^{n},can)$, define

$$V_{PC}(n,\varphi) = \underset{\mu\in \Mob(\SS^{n})}{\sup} \Vol(\mathcal{O}, (\mu\circ \varphi)^{*}(can)).$$
If there exists a PC map $\varphi: |\mathcal{O}| \to \SS^{n}$, then we also define
$$V_{PC}(n, \mathcal{O}) = \underset{\varphi:|\mathcal{O}|\to \SS^{n}  PC}{\inf} V_{PC}(n,\varphi).$$

$V_{PC}(n,\mathcal{O})$ is denoted the ($n$-dimensional) {\it piecewise conformal volume} of $\mathcal{O}$.
\end{definition}

{\bf Remark:} It seems likely that our definition of piecewise conformal volume coincides with conformal volume of Li-Yau for manifolds,
but we have not checked this (it would suffice to show that a PC map can be approximated
by conformal maps).

If there exists a piecewise isometric map $\varphi: (|\mathcal{O}|,g) \to \EE^{n}$ for some $n$,
then clearly $V_{PC}(n,\OO)$ is well-defined, since $\EE^{n}$ has a conformal embedding
into $\SS^{n}$. Let $\S\subset \OO$ be the singular locus of $\OO$. By the Nash embedding Theorem \cite{GromovRohlin70},
there is an isometry $\rho: \OO\backslash \S \to \EE^{N}$, for $N$ large enough. Then $\rho$
extends continuously to an isometric map $\overline{\rho}: \OO \to \EE^{N}$.
Thus, $V_{PC}(N, \OO)$ is well-defined.
We record some basic facts about conformal volume, generalizing some of the facts in \cite[\textsection 1]{LiYau82}
which carry over to our notion of conformal volume for orbifolds.

{\bf Fact 1:} If $\mathcal{O}$ admits a PC map $\psi$ of geometric  degree $d$ onto another orbifold $\mathcal{P}$, then
$$V_{PC}(n,\mathcal{O}) \leq |d| V_{PC}(n,\mathcal{P}).$$ This follows because for any PC
map $\varphi: \mathcal{P} \to \SS^{n}$, the composition $\varphi\circ \psi: \OO \to \SS^{n}$
satisfies $V_{PC}(n,\varphi\circ \psi) \leq d V_{PC}(n,\varphi)$.

{\bf Fact 2:} Since the embedding $\SS^{n}\hookrightarrow \SS^{n+1}$ is an isometry, it's clear that $V_{PC}(n,\mathcal{O})\geq V_{PC}(n+1,\mathcal{O})$.
Define the piecewise conformal volume $V_{PC}(\mathcal{O})=\underset{n\to \infty}{\lim} V_{PC}(n,\mathcal{O})$.

{\bf Fact 3:} If $\mathcal{O}$ is of dimension $m$, and $\varphi: |\mathcal{O}| \to \SS^{n}$ is a PC map, then
$$V_{PC}(n,\varphi)\geq V_{PC}(n,\SS^{m})=\Vol(\SS^{m}).$$
The same argument as in Li-Yau works here: ``blow up'' about a smooth  manifold point of $\varphi(|\mathcal{O}|)$
so that  the image converges in the Hausdorff topology to a geodesic sphere of dimension $m$. More precisely,
let $x\in \OO$ be a manifold point of $\OO$ such that the differential $d\varphi_{x}$ has rank $m$. Let $\mu_{t}:\RR^{n+1}\to \RR^{n+1}$
be a M\"obius transformation preserving $\SS^{n}$, fixing $\pm\varphi(x)$, and translating the origin $0\in \RR^{n+1}$ to $-t \varphi(x)$, $0\leq t <1$.
Let $U\subset \OO $ be a small neighborhood of $x$ such that $rk(d\varphi_{y})=m$, for all $y\in U$. Then
$\mu_{t}(\varphi(U))$ converges to the geodesic sphere of dimension $m$ going through $\pm \varphi(x)$, and
tangent to $d\varphi_{x}(T_{x}\OO)$ as $t\to 1$. Thus, we see that $V_{PC}(n,\varphi) \geq \Vol(\SS^{m})$.

{\bf Fact 4:} If $\mathcal{O}$ is an embedded suborbifold of the orbifold $\mathcal{P}$, and $\varphi: |\mathcal{P}|\to \SS^{n}$
is  PC, then $V_{PC}(n, \varphi) \geq V_{PC}(n,\varphi|_{|\mathcal{O}|})$.
Thus, $V_{PC}(n,\mathcal{P}) \geq V_{PC}(n, \mathcal{O})$.

For example, suppose $\i: P\subset \SS^{n}$ is an $n$-dimensional submanifold, with the conformal structure induced from the
embedding. Then $V_{PC}(n,\iota) \leq V_{PC}(\SS^{n})$, by fact 4.
But by fact 3, $V_{PC}(n,\iota) \geq V_{PC}(\SS^{n})$, and therefore
we have shown that $V_{PC}(n, \iota) = V_{PC}(P)=V_{PC}(\SS^{n}) =\Vol(\SS^{n})$.

{\bf Fact 5:} If $\OO= \OO_{1} \cup \OO_{2}$, and $\varphi:\OO\to \SS^{n}$,
then $$V_{PC} (\mathcal{O}, \varphi) \le V_{PC} (\mathcal{O}_1, \varphi|_{\mathcal{O}_1}) + V_{PC} (\mathcal{O}_2, \varphi|_{\mathcal{O}_2} ).$$

\section{Finite subgroups of $\O(n+1)$}

In this section, we give a bound on the conformal volume of $\SS^{n}/\G$, where $\G<\O(n+1)$
is a finite subgroup. For $\G < \O(n+1)$ finite, define the orbifold $M^n_\G = \mathbb{S}^n / \G.$

\begin{theorem} \label{ellipticconformalvolume}
There is a function $C(n)$ such that if $\G< \O(n+1)$, then
$$V_{PC}(M^{n}_{\G}) < C(n) |\G|^{2n}.$$
\end{theorem}
\begin{proof}
First, we need a lemma.
\begin{lemma} \label{slice}
If we have $k$ hyperplanes in $\RR^{n}$ going through the origin, where $k>0, n>0$,
then they separate $\RR^{n}$ into at most $f(k,n)=2\sum_{i=0}^{n-1} \binom{k-1}{i} = O(k^{n-1})$ regions.
\end{lemma}
\begin{proof}
We may assume the hyperplanes are in general position.
The proof is by lexicographic induction on $(k,n)$. Clearly, $f(k,1) = 2$,  $f(1,n)=2$.
Consider $k$ hyperplanes $\{P_{1}, ..., P_{k}\}$, $P_{i}\subset \RR^{n}$, $k>1, n>1$.
Then the hyperplanes $\{P_{1},..., P_{k-1}\}$ separate $\RR^{n}$ into at most
$f(k-1,n)$ regions, by the inductive hypothesis. Also,
$\{P_{1}\cap P_{k}, ..., P_{k-1}\cap P_{k}\}$ is a collection of $k-1$ hyperplanes in $P_{k}\cong \RR^{n-1}$,
and thus they divide $P_{k}$ into at most $f(k-1, n-1)$ regions,
by the inductive hypothesis. Each region of $\RR^{n}\backslash (P_{1}\cup ... \cup P_{k-1})$
which meets $P_{k}$ will be divided into two, and there are
at most $f(k-1,n-1)$ such regions. The regions which don't
meet $P_{k}$ will remain unchanged. Thus,
$|\RR^{n}\backslash (P_{1}\cup ...\cup P_{k})| \leq f(k-1,n) + f(k-1,n-1)=f(k,n)$.
This follows since
$$f(k-1,n) +f(k-1,n-1)= 2 \sum_{i=0}^{n-1} \binom{k-2}{i} + 2\sum_{i=0}^{n-2} \binom{k-2}{i}$$
$$=2 \sum_{i=0}^{n-1} \binom{k-2}{i} + \binom{k-2}{i-1} = 2\sum_{i=0}^{n-1} \binom{k-1}{i} = f(k,n).$$

\end{proof}

We will actually find a piecewise isometric map $f:\SS^{n}/\G \to \SS^{N}$ for some $N$, in which
the number of isometrically embedded pieces is bounded by a polynomial in $|\G|=m$. This will immediately
imply the result, using facts 3 through 5 about conformal volume.

Define the orbi-covering $\pi_\G: \mathbb{S}^n \longrightarrow M^n_\G$.
Let $\G=\{ \g_1, \g_2, \ldots, \g_m \}$.  
For each element $\gamma \in \Gamma$,  right multiplication by $\g$ induces a permutation of $\{\g_{1},\g_{2}, \ldots, \g_{m} \}$, 
which we identify with a permutation of $\{1,2,\ldots, m\}$.  
Notice this does not give an action of $\Gamma$, since composing right multiplications 
reverses the order of group multiplication. Use these permutations to define, for each $\g\in \G$,
a permutation  of the coordinates of $\prod_{i=1}^m \overline{\mathbb{B}}^{n+1}$ given by
$$\g \cdot (x_1, \ldots, x_m ) =  (x_{\g \cdot 1}, \ldots, x_{\g \cdot m}).$$

Define $N = m (n+1)-1$ and
\begin{eqnarray*}
\Psi : \SS^{n}  & \longrightarrow & \mathbb{S}^N \subset \prod_{i=1}^m \overline{\mathbb{B}}^{n+1}   \\
x  & \longmapsto &  m^{-\frac12}( \g_1 x,  \g_2 x, \ldots , \g_m x).
\end{eqnarray*}
$\SS^{n}$ is taken by $\Psi$ to a totally geodesic submanifold of $\SS^N$,   since it is (up to rotation in $\O(N+1)$) the
restriction of the diagonal map $\RR^{n+1} \to \RR^{N+1}$, rescaled to be an isometry.
For each element $\g\in \G$
the map $\Psi$ is $\g$-equivariant, because  
$$\Psi( \g\cdot x) =m^{-\frac12} (\g_{1}\cdot \g \cdot x, \ldots, \g_{m} \cdot \g\cdot x)=m^{-\frac12}(\g_{\g\cdot 1}x,\ldots,\g_{\g\cdot m}x)= \g\cdot \Psi(x) $$ 
for any  $x\in \SS^{n}$.

Let $\Sigma_{m}$ be the permutation group on $m$ elements.
The map $\Psi$ therefore descends to a map
$$ \Upsilon: M^{n}_\G \longrightarrow \SS^N / \S_{m} = Y.$$

The action of $\S_{m}$ on $\SS^N$ induces an embedding $\S_{m}< \O(N+1)$.  By permuting the coordinates of $\RR^{N+1}$, the symmetric group $\Sigma = \S_{N+1}$ acts isometrically as a reflection group on $\SS^N$, such that $\S_{m}<\S$.  Similar to what we did above, define the orbi-covering $\pi_\S: Y \longrightarrow M^N_\S$ and an isometric embedding $\iota_\S: |M^N_\S| \longrightarrow \SS^N$, which exists since we may identify $|M^{N}_{\S}|$ with
the Coxeter polyhedron for $\S$.
Let $\{ S_k\}_{k \in K}$ be the collection of $\binom{N+1}{2}$ hyperspheres given by the fixed point sets of the reflections in $\S$.

\emph{The goal is now to estimate $V_{PC} ( M^{n}_\Gamma, ( \iota_\Sigma \circ \pi_\Sigma \circ \Upsilon)) $ from above.}

Consider how $\Psi (\mathbb{S}^n)$ is cut up by the collection $\{ S_k \}_{k \in K}$.
It is cut up into at most $f\left(\binom{N+1}{2},n+1\right)$ convex sub-polyhedra $\{ C_\s, \s\in \S\}$
by the planes $\{\Psi(\SS^{n}) \cap S_k \}_{k \in K}$ (For most $\s\in \S$, $C_{\s}=\emptyset$).  Then each $C_\s$ is isometrically and totally geodesically projected to $\SS^{N}$ via the map $\i_{\S}\circ\pi_{\S}\circ\Upsilon$.  Postcomposition by a conformal automorphism of $\mathbb{S}^N$ can then make the volume of the image of $C_\s$ no more than the volume of a unit $n$-sphere. By the combination of facts
3 through 5 of section \ref{conformal volume}, we can conclude that
$$V_{PC} (M^{n}_\Gamma, ( \iota_\Sigma \circ \pi_\Sigma \circ \Upsilon )) \le f\left(\binom{m(n+1)}{2} ,n+1\right) \cdot \text{Vol}(\mathbb{S}^n)\leq C(n) |\G|^{2n},$$
for some constant $C(n)$.
\end{proof}

{\it Remark:} It turns out that a finite subgroup $\G < \O(n+1)$ does not necessarily embed in a finite reflection
subgroup of $\O(n+1)$. Thus, the strategy of \cite{Agol05} does not generalize to higher dimensions. It is an interesting
problem to try to get better estimates on $V_{PC}(M^{n}_{\G})$. The method we have used to estimate
this quantity is probably far from optimal. In fact, in the above theorem we have constructed maps which potentially have
very high multiplicity, which is exactly the sort of thing one  wants to avoid in estimating conformal volume.

\section{Eigenvalue bounds}

In this section, we observe that the argument of \cite[Thm. 2.2]{ElSoufiIlias86} generalizes to our
context of orbifolds (their theorem sharpens Corollary 3 from section 2 of \cite{LiYau82}). If $(\OO,g)$
is a Riemannian orbifold, then $\l_{1}(\OO,g)$
is the first non-zero eigenvalue of the operator $-\tilde{\D}_{g}$ on $\mathcal{O}$, where 
$\tilde{\D}_{g}$ is the unique self-adjoint extension of the Laplacian on $\mathcal{O}$. 
We will be using the variational characterization of $\l_{1}(\OO,g)$ in
terms of the Rayleigh-Ritz quotient.

\begin{lemma}
Let $\OO$ be an orbifold of finite volume. Let $f\in W^{1,2}(\OO)$
be an $L^{2}$ function with $L^{2}$ distributional derivative, and such that $\int_{\OO} f dv_{g}=0$. Then 
$$\l_{1}(\OO,g)\int_{\OO}f^{2}dv_{g}\leq \int_{\OO}|\nabla f|^{2} dv_{g}.$$
\end{lemma}

\begin{proof}
Let $A$ be a self-adjoint operator on a Hilbert space $\mathcal{H}$, such that $A$ is semi-bounded
below. This means that there exists a constant $k$, such that for all $\psi\in \mathcal{H}$,
$\langle \psi, A\psi \rangle \geq k \langle \psi, \psi \rangle$. The Laplacian $\D$ with domain
restricted to compactly supported, $C^{\infty}$ functions on $\OO$ is essentially
self-adjoint (\cite[Thm. 1.9]{EGM98}), with self-adjoint extension $\tilde{\D}$. 
 In the case of the self-adjoint
operator $-\tilde{\D}$, we have $-\tilde{\D}$ is semi-bounded below, with constant $k=0$. This
follows from the fact that $C_{c}^{\infty}(\OO)$ is dense in 
$Domain(\tilde{\D})$, and for $f\in C_{c}^{\infty}(\OO)$ integration by parts 
gives $\int_{\OO} -f \D f dv_{g}= \int_{\OO} |\nabla f|^{2} dv_{g}$, which is non-negative. 
We have $\l_{0}(-\tilde{\D})=0$, with constant eigenfunction. Therefore, the
Rayleigh-Ritz principle for semi-bounded below operators \cite[XIII]{ReedSimon78} gives 
$$\l_{1}(\OO,g) = \inf \left\{ \frac{\int_{\OO} -f\tilde{\D} f dv_{g}}{\int_{\OO} f^{2} dv_{g}}\, |\, f\in Domain(\tilde{\D}),  \int_{\OO} f dv_{g}=0\right\}.$$ 
Since $C_{c}^{\infty}(\OO)$ is dense in $Domain(\tilde{\D})$, we may take the above
infimum over functions $f\in C_{c}^{\infty}(\OO)$. Then integrating by parts, 
we have
$$\l_{1}(\OO,g) = \inf \left\{ \frac{\int_{\OO} |\nabla f|^{2} dv_{g}}{\int_{\OO} f^{2} dv_{g}}\, |\, f\in C_{c}^{\infty}(\OO),  \int_{\OO} f dv_{g}=0\right\}.$$ 
Then we use again the fact that $C_{c}^{\infty}(\OO)$ is dense in the Sobolev space
$W^{1,2}(\OO)$ to see that we may take the above infimum over functions in
$W^{1,2}(\OO)$. 
\end{proof}

\begin{theorem} (see \cite[Thm. 2.2]{ElSoufiIlias86}) \label{ev}
Let $(\mathcal{O},g)$ be a Riemannian orbifold of dimension $m$ of finite volume. 
If $\varphi:|\mathcal{O}|\to \SS^{n}$ is a PC map,
$$\l_{1}(\mathcal{O},g)\, \Vol(\mathcal{O},g)^{\frac2m} \leq m V_{PC}(n,\varphi)^{\frac2m}.$$
\end{theorem}
\begin{proof}
Let ${\bf X} = (X_{1}, ..., X_{n+1})$ be the coordinate
functions on $\RR^{n+1}$. Then $\sum_{i=1}^{n+1} X_{i}^{2}= 1$, restricted to the unit sphere $\SS^{n}\subset \RR^{n+1}$.

\begin{lemma}
There exists $\mu \in \Mob(\SS^{n})$ such that
$\int_{\mathcal{O}} {\bf X}\circ \mu \circ \varphi\ dv_{g}={\bf 0}$.
\end{lemma}

\begin{proof}
Let ${\bf x} \in \SS^{n}$. For $0\leq t<1$, $t {\bf x}\in \BB^{n+1}$ (the open unit ball),
let $\mu_{t{\bf x}} \in \Mob(\SS^{n}) = \Isom(\BB^{n+1})$
be the hyperbolic translation along the ray $\RR{\bf x}$ taking ${\bf 0}$
to $t{\bf x}$ (thus, $\mu_{0 {\bf  x}}=\mu_{{\bf 0}}=Id$).
Let $H(t{\bf x}) = \frac{1}{\Vol(\mathcal{O},g)} \int_{\mathcal{O}} {\bf X}\circ \mu_{t{\bf x}}\circ \varphi\ dv_{g}$.
We may think of $H$ as defining a continuous function $H:\BB^{n+1}\to \BB^{n+1}$,
which gives the center of mass of the measure
coming from $\varphi_{*} dv_{g}$, where
we take the point $-t{\bf x}$ to the origin of the sphere
$\SS^{n}=\partial \BB^{n+1}$ by the conformal map $\mu_{t{\bf x}}$. Let $x_{i}\in \BB^{n+1}$
be a sequence of points approaching $x_{\infty}\in \SS^{n}$. We want to
show that $H$ extends to a map fixing $\partial \BB^{n+1}$. To see this, let
$U_{\delta}=\SS^{n}\backslash B_{\delta}(-x_{\infty})$ be a neighborhood of $x_{\infty}$.
Then the fact that $\varphi$ is PC implies that for any $\e>0$, there exists $\delta>0$,
such that $\Vol(\varphi^{-1}(U_{\delta})) > \Vol(\OO)(1-\e)$. Then we may find $N>0$, such
that $\mu_{x_{i}}(U_{\delta}) \subset B_{\e}(x_{\infty})$, for $i>N$. Thus, we may see that
$H(x_{i}) \in B_{\e}(x_{\infty})$. So $\underset{i\to \infty}{\lim} H(x_{i}) = x_{\infty}$.
So we see that $H$ extends to a continuous
function $\overline{H}:\overline{\BB}^{n+1}\to \overline{\BB}^{n+1}$  such that
$\overline{H}|_{\SS^{n}}=Id|_{\SS^{n}}$, and $\overline{H}|_{\BB^{n+1}}=H$.  Thus, $H$ is onto, so there exists ${\bf y} \in \BB^{n+1}$
such that $H({\bf y}) = {\bf 0}$, and we take $\mu=\mu_{{\bf y}}$. \end{proof}

Now, replace $\varphi$ with $\mu\circ \varphi$, noting that this
is still a PC map. Then $X_{i}\circ \varphi$ may be
used as test functions in the Rayleigh quotient, since they
are $L^{2}$ orthogonal to the
constant function, and lie in the Sobolev space $W^{1,2}(\OO)$ since they
are piecewise smooth (they are differentiable {\it a.e.}, with (distributional) derivative lying in 
$L^{2}(\OO)$, which follows from equations (2) and (3) below, if $\vol(\OO,\varphi^{*} can)<\infty$, which
we may assume if $V_{PC}(\OO, \SS^{n})< \infty$). Thus,
\begin{equation}
\l_{1}(\mathcal{O}) \int_{\mathcal{O}} |X_{i}\circ \varphi|^{2} dv_{g}\leq \int_{\mathcal{O}} |\nabla X_{i}\circ \varphi |^{2} dv_{g}.
\end{equation}
Summing, we see that
\begin{equation*}
\l_{1}(\mathcal{O}) \Vol(\mathcal{O},g) = \l_{1}(\mathcal{O}) \int_{\mathcal{O}} \sum_{i=1}^{n+1}   |X_{i}\circ \varphi|^{2} dv_{g} \leq m \int_{\mathcal{O}} \frac1m \sum_{i=1}^{n+1} |\nabla X_{i}\circ \varphi |^{2} dv_{g}
\end{equation*}
\begin{equation}\label{eqn1}
 \leq m \left(\int_{\mathcal{O}} \left( \frac1m \sum_{i=1}^{n+1} |\nabla X_{i}\circ \varphi|^{2}\right)^{\frac{m}{2}} dv_{g}\right)^{\frac2m} \Vol(\mathcal{O},g)^{1-\frac2m} ,
 \end{equation}
where the last inequality is H\"older's inequality.

Take a point $x\in \OO$  at which $\varphi$ is differentiable and conformal (By the definition of a PC map,
almost every point of $\OO$ satisfies these two conditions). After possibly performing an orthogonal change of coordinates at $x \in \mathcal{O}$ and on $\mathbb{R}^{n+1}$, we may assume without loss of generality that the derivative $d\varphi(x): T_{x}(\OO) \to T_{\varphi(x)}(\RR^{n+1})$ may be written
as
$$d \varphi = \left[\begin{array}{ccc}k & 0 & 0 \\0 & \ddots & 0 \\0 & 0 & k \\0 & 0 & 0 \\0 & 0 & 0\end{array}\right],$$
where $k$ is the dilatation factor of $\varphi$ at $x$. From the fact that $\tr( A^{t}A)$ is invariant under pre- or post-composition with an orthogonal matrix,
one sees that
$$  \frac1m \sum_{i=1}^{n+1} |\nabla X_{i}\circ \varphi |^{2} =\frac1m \tr(d\varphi(x)^{t} d\varphi(x) ) = k^{2}.$$
Since $\varphi^{*}(can)(x) = k^{m} dv_{g} (x)$, where $can$ is the volume form on $\SS^{n}$, we see that
$$\left(\frac1m \sum_{i=1}^{n+1} |\nabla X_{i}\circ \varphi |^{2}\right)^{m/2} dv_{g} = \varphi^{*}(can)\,  \text{a.e.} $$
This yields the inequality
\begin{equation} \label{eqn2}
\int_{\mathcal{O}} \left( \frac1m \sum_{i=1}^{n+1} |\nabla X_{i} \circ \varphi |^{2} \right) ^{\frac{m}{2}} dv_{g}  \leq \Vol(\mathcal{O},\varphi^{*}can).
\end{equation}
Finally, we obtain the desired inequality by combining inequalities \eqref{eqn1} and \eqref{eqn2}, and
taking an infimum over all PC functions $\varphi$.
\end{proof}

\section{Congruence arithmetic hyperbolic orbifolds}
This section collects some properties of arithmetic hyperbolic orbifolds which
will be needed in the next section. Given a number field $k$, let $\cO$ denote
its ring of integers, $\disc_k$ the absolute value of the discriminant of $k$
and $\deg(k) = [k:\QQ]$.
\medskip

The group $H = \Isom(\HH^{n})$ can be identified with $\O_0(n,1)$ --- the
subgroup of the orthogonal group $\O(n,1)$ which preserves the upper-half
space. We now define arithmetic and congruence subgroups of $H$.

\def\ia{{\mathfrak a}}

\begin{definition}
Let $k$ be a totally real number field, $f$ a quadratic form of signature
$(n,1)$ defined over $k$ such that for every non-identity embedding
$\sigma :k\to\RR$ the form $f^\sigma$ is positive definite. Then the group
$\Gamma = \O_0(f,\cO)$ of the integral automorphisms of $f$ is a discrete
subgroup of $H$ via the inertia theorem. Such groups $\Gamma$ and subgroups of $H$ which are
commensurable to them are called {\it arithmetic subgroups of the simplest
type}. The field $k$ is called a {\it field of definition} of $\Gamma$ (and
subgroups commensurable to $\G$ in $H$). An arithmetic subgroup $\Gamma$ is called a {\it
congruence subgroup} if there is a nonzero ideal $\ia\subset\cO$ such that
$\Gamma\supset \O_0(f,\ia)$, where $\O_0(f,\ia) = \{ g\in\O_0(f,\cO) \mid
g\equiv Id\ (\mod\ \ia)\}$, the {\it principal congruence} subgroup of
$\O_0(f,\cO)$ of level $\ia$ (See \cite{PR} for more on arithmetic groups).
We will also apply this terminology to the corresponding quotient orbifolds.
\end{definition}

For the hyperbolic spaces of even dimension all arithmetically defined
subgroups are arithmetic subgroups of the simplest type. For
odd $n$ there is another family of arithmetic subgroups given as the groups of
units of appropriate Hermitian forms over quaternion algebras. Moreover, if
$n=7$ there is also the third type of arithmetic subgroups of $H$ which are
associated to the Cayley algebra. This description follows from Tits'
classification of semisimple algebraic groups \cite{Tits66}. For our purpose it
will be sufficient to consider only arithmetic subgroups of the simplest type, due to the following result of Vinberg:
\begin{lemma}\cite[Lemma~7]{Vinberg67} 
 If $\Gamma$ contains an arithmetic subgroup
generated by reflections then it is defined by a quadratic form.
\end{lemma}

The following lemma is a special case of Proposition~3.3 in \cite{Bel} which
was proved as an application and extension of the method of Borel-Prasad
\cite{BorelPrasad}.
\begin{lemma} \label{volbydisc}
Given $n>3$, there exist positive constants $c_{n} , \delta_{n}$ such that if a
maximal arithmetic subgroup $\G < \Isom(\HH^{n})$ is defined over the field
$k$, then $\Vol(\HH^{n}/\Gamma) \ge c_{n}\disc_k^{\delta_{n}}$.
\end{lemma}
\begin{proof}
It was shown in \cite[Sect.~3.3]{Bel} that for almost all admissible fields
$k$, if $\Gamma$ is an arithmetic subgroup of $H (=\Isom(\HH^{n}))$ defined
over $k$, then
% \marginpar{I've made a slight modification in \cite{Bel}, so this inequality can be found in sect. 3.3 now.}
$$ \mu(H/\Gamma) \ge \frac12\disc_k^{\delta_n},$$
where $\mu$ is a Haar measure on $H$ normalized in a certain way and $\delta_n$
is a positive constant which depends only on $n$ and can be computed
explicitly. There is a constant $c = c(n)$ such that $\Vol(\HH^{n}/\Gamma) =
c\mu(H/\Gamma)$ (in fact, in our case $c$ is equal to the volume of the unit
sphere in $\RR^n$). Thus for all but finitely many $k$,
$\Vol(\HH^{n}/\Gamma)\ge \frac12c\disc_k^{\delta_n}$. We can extend this
inequality to the remaining finite collection of fields by decreasing if
necessary the constant $c$.
\end{proof}

The next fact which we recall here may have been known for a long time. More recently, related questions
in a different context were thoroughly studied by I.~Horozov in his PhD thesis \cite{Horozov, Horozov05}.
Still we were not able to find a proof of the statement in the literature, so we will present an argument.
\begin{lemma} \label{cyclicvolume}
There exists a constant $e_{n}$ such that if an arithmetic subgroup $\Gamma <
\Isom(\HH^{n})$ is defined by a quadratic form defined over $k$, and $C
< \Gamma$ is a finite cyclic subgroup, then $|C| \le e_{n} \deg(k)^{2n+2}$.
\end{lemma}

\begin{proof}
%We can assume that $\Gamma$ is a maximal arithmetic subgroup of $H$ defined
%over $k$. There exists a $k$-group $G$ and a $k$-embedding
%$G\hookrightarrow\GL(m)$ such that $\Gamma = N_H(\Gamma\cap G(k))$. Moreover,
%$\Gamma/\Gamma\cap G(k)$ is an elementary abelian group and its exponent
%divides the order of the center of the simply connected cover of $G$. These
%facts are well known and not hard to check.
%
%Let $\Lambda = \Gamma\cap G(k)$. In our case $\Gamma/\Lambda$ is an elementary
%abelian $2$-group. Thus the torsion of order $>2$ appears already in
%$\Lambda\subset G(k)$. We come to the question about the torsion in
%$G(k)\subset\GL(m,k)$ and as $\Gamma$ is defined by a quadratic form $m= n+1$.
%\marginpar{before I claimed that we can assume $\Lambda\subset\GL(m,\cO)$ but
%this is not always the case! fortunately, this assumption is not needed.}

Let $\Gamma$ be a subgroup of $H$ commensurable with the group of integral
points of $G(k) = \O_0(f,k)$. It is known that if the center of $H$ is trivial
(which is indeed the case for $H = \Isom(\HH^{n})$) then $\Gamma$ is contained
in $G(k)$ \cite[Prop. 1.2]{ChernousovRyzhkov97}. We come to a question about the torsion in
$G(k)\subset\O_0(f,k)\subset\GL(n+1,k)$.

Let $A\in \GL(m,k)$, $m = n+1$ be a torsion element of order $t\ge2$. Let
$\lambda_1, \ldots,\lambda_m$ be the eigenvalues of $A$. As $A$ is a matrix
over $k$, its eigenvalues split into groups of conjugates under the action of
$Gal(\bar{k}/k)$. Now, $A^t= Id$ implies that the eigenvalues are roots of
unity. Let $t_1, \ldots, t_l$ be their orders, so $t={\rm lcm}(t_1, \ldots,
t_l)$. If $\lambda$ is an eigenvalue then all its $Gal(\bar{k}/k)$-conjugates
are also eigenvalues of $A$, which implies
$$
\phi_k(t_1)+\cdots+\phi_k(t_l) \le m,
$$
where $\phi_k(t)$ denotes a generalized Euler $\phi$-function, which can be defined as the
degree over $k$ of the cyclotomic extension $k(\mu_t)$, where $\mu_t$ denotes a primitive
$t$th-root of unity.

It is clear that the following inequalities are satisfied for $\phi_k(t)$ and the Euler
$\phi$-function:
$$
\phi(t)/\deg(k)\le \phi_k(t) \le \phi(t).
$$
To simplify the notation let $\deg(k)=d$. We have
$$
\phi(t_i)\le md.
$$
This implies
$$
\phi(t) \le \phi(t_1)\cdots\phi(t_l) \le (md)^m
$$
(the first inequality employs the fact that $t = {\rm lcm}(t_1, \ldots, t_l)$).

Now using the well known inequality $\phi(t)\ge \sqrt{t}/2$, we obtain
$$
t\le 4 (md)^{2m} \le e_{n}d^{2n+2}.
$$
% In order to come back to the torsion in a maximal arithmetic subgroup $\Gamma <
% H$ we may need to multiply $t$ by $2$ which can be captured by slightly
% enlarging the constant.
\end{proof}

\begin{cor} \label{finitesubgroup}
Under the hypothesis of the above lemma, there exists a constant $m_{n}$ such that
if $F< \G$ is a finite subgroup, then $|F|\le m_{n} \deg(k)^{n(n+1)}$.
\end{cor}
\begin{proof}
%\marginpar{Rem: what we use here is a quantitative version of a theorem of
%Jordan proved by Frobenius but I am not sure about the reference. We can also
%leave it as it is.}
By Margulis' Lemma, there is a constant $M_{n}$ such that
if $F< \O(n)$, then there is an abelian subgroup $A< F$, such that $[F:A] \leq
M_{n}$ (see {\it e.g.} \cite[Cor. 4.1.11]{Thurston97}). We may find common
complex eigenspaces of the elements of $A$: $U_{1}, \ldots, U_{k}$, and real
eigenspaces $V_{1},\ldots,V_{l}$, where $k/2+l\leq n$, such that $A$ acts on
$U_{i}$ as a cyclic group $A_{i}$, and $A$ acts on $V_{i}$ as $\pm 1$. We may
embed $A$ in $\prod_{i=1}^{k} A_{i} \times \left(\ZZ/2\ZZ\right)^{l}$, acting on
$\prod_{i=1}^{k} U_{1} \times \prod_{j=1}^{l} V_{j}$.  Thus $|A| \leq
2^{l}\prod_{i=1}^{k} |A_{i}|$. By the previous lemma, $|A_{i}| \leq
e_{n}\deg(k)^{2n+2}$, since a generator of $A_{i}$ is the projection of an
element of $A$. Thus, $|F| \leq M_{n}|A| \leq m_{n} \deg(k)^{n(n+1)}$.
\end{proof}

\begin{lemma} \label{discriminant}
There exists $c>1$, such that for a number field $k$, $\disc_k > c^{\deg(k)}$.
\end{lemma}

This follows from Minkowski's theorem (see {\it e.g.} \cite[V.4]{LangANT}), however
if one wants to obtain a ``good value'' for the constant $c$, the question becomes much more
complicated. There has been a considerable amount of research in this direction. We refer to
\cite{Od90} for a survey of the results available so far.

\begin{lemma}
Maximal arithmetic subgroups of the simplest type are congruence.
\end{lemma}
\begin{proof}
Let $\Gamma$ be a maximal arithmetic subgroup which is commensurable with
$G(\cO)=\O_0(f,\cO)$ for some quadratic form $f$ defined over $k$.
There exists an $\cO$-lattice $L$ in $V=k^{n+1}$ such that $\Gamma\cap G(k)\
(=\Gamma)= G^L = \{g\in G(k) \mid g(L) = L\}$ (see \cite[Prop.~4.2]{PR} for
$k = \QQ$, the general case is entirely similar as we do not claim that $L$ is
a free $\cO$-lattice). Let also $M$ be the standard lattice in $V$ so that $G^M
= G(\cO)$. We claim that there exists a principal congruence subgroup of $G^M$
which is contained in $G^L$. The argument is the same as in
\cite[Lemma~4.1.1]{PR}: Let $\phi: G\to G$ be a $k$-isomorphism defined by the
change of basis of $V$ from the standard basis (of $M$) to a $k$-basis
contained in $L$. It follows that if $\Delta$ is a principal congruence
subgroup of $G^M$ with respect to an ideal $\ia\subset\cO$ such that the
coefficients of $\ia (\phi-Id)$ in the standard basis are integral, then
$\Delta\subset G^L$.
\end{proof}

{\bf Remark:} The results which were discussed in this section so far are
valid for arbitrary arithmetic subgroups, the proofs for the general
case are similar to the case of the simplest type which we considered in
detail.

The following theorem was stated in the case of arithmetic groups defined over
$\QQ$, but is true in general. It is a translation of results of
Jacquet-Langlands \cite{JL70} and Gelbart-Jacquet \cite{GJ78} in the theory of
automorphic forms into a statement about eigenvalues of congruence arithmetic
groups, and it generalizes a famous theorem of Selberg for congruence subgroups
of $\PSL(2,\ZZ)$ \cite{Selberg65}. 
\begin{theorem} [Vigneras \cite{Vigneras83}]
Let $\OO = \HH^{2}/\G$, where $\G$ is an arithmetic congruence Fuchsian
subgroup of $\Isom(\HH^{2})$. Then $\l_{1}(\OO)\geq \frac{3}{16}$.
\end{theorem}

In the above theorem, it is conjectured that $\l_{1}(\OO) \geq \frac14$, which
generalizes a conjecture of Selberg about congruence subgroups of
$\PSL(2,\ZZ)$. The following theorem of Burger and Sarnak extends this result
to higher dimensions, by inducing from the 2-dimensional case.

\begin{theorem} \cite[Cor. 1.3(a)]{BurgerSarnak91} \label{BurgerSarnak}
Let $\OO=\HH^{n}/\G$, where $\G$ is an arithmetic congruence subgroup defined
by a quadratic form, $n\geq 3$ (or if $n=3$, $\G$ need only be arithmetic
congruence). Then $\l_{1}(\OO) \geq \frac{2n-3}{4}$.
\end{theorem}
Under the hypotheses of Theorem \ref{BurgerSarnak}, it is conjectured that
$\l_{1}(\OO)\geq n-2$, which is a special case of the {\it generalized
Ramanujan conjecture} \cite[Cor. 1.3(b)]{BurgerSarnak91}.

\section{Finiteness of arithmetic hyperbolic maximal reflection groups}

We put together the results in the previous sections to prove our main theorem.

\begin{theorem}
There are only finitely many conjugacy classes of arithmetic maximal reflection groups in
${\rm Isom}(\HH^{n})$.
\end{theorem}
\begin{proof}

Suppose that $\G$ is an arithmetic maximal reflection group. Then there
is no group $\G' < \Isom(\HH^{n})$ such that $\G <\G'$ and $\G'$ is generated by reflections.
There  does exist a maximal lattice $\G_{0}<\Isom(\HH^{n})$ satisfying $\G\leq \G_{0}$.
The group $\G$ is generated by reflections in a finite volume polyhedron $P$, which forms
a fundamental domain for $\G$.

\begin{lemma} [Vinberg \cite{Vinberg67}] \label{vinb2}
$\G$ is a normal subgroup of $\G_{0}$. Moreover, there is a finite subgroup $\T<\G_{0}$
such that $\T \to  \G_{0}/\G$ is an isomorphism, and $\T$ is the group of symmetries of
 the polyhedron $P$.
\end{lemma}
\begin{proof}
This  follows from the fact that the set of reflections in $\G_{0}$
is conjugacy invariant, and therefore the subgroup of $\G_{0}$ generated by reflections in $\G_{0}$ is normal
in $\G_{0}$. Since $\G$ is a maximal reflection group, this subgroup of $\G_{0}$ must be $\G$.

Suppose there is an element
$\g\in \G_{0}$ such that $\interior(P)\cap \g(\interior(P)) \neq \emptyset$ and $\g(P)\neq P$.
Then there is a geodesic plane $V$ containing a face of $P$ such
that $\g(V)\cap \interior(P) \neq \emptyset$. The reflection $r_{V}\in \G$ in the plane $V$
is conjugate to a reflection $r_{\g(V)}= \g r_{V} \g^{-1}$, which is not
in $\G$ since $r_{\g(V)} (\interior(P))\cap \interior(P) \neq \emptyset$. However
$\G\unlhd \G_{0}$ implies $r_{\g(V)} \in \G$, a contradiction.

Let $\T$ be the subgroup of
$\G_{0}$ such that $\T(P)=P$.  $\T$ is finite because $P$ is finite volume and
has finitely many faces.

Consider the composition $\T \hookrightarrow \G_{0} \rightarrow \G_{0}/\G$. $\G\cap \T$ is trivial because
$P$ is a fundamental domain for $\G$. The composition map is therefore injective. To see surjectivity,
pick a $\g_{0} \in \G_{0}$ and find a $\g\in \G$ such that $\g_{0}^{-1} (\interior(P)) \cap \g(\interior(P)) \neq \emptyset$.
Then $\g_{0}\g(P)=P$, so $\g_{0} \g = \th \in \T$. Therefore $\g_{0}\G = \th \G$. This shows that the
composition $\T \rightarrow \G_{0}/\G$ is surjective, and an isomorphism.
\end{proof}

Let $\mathcal{O}=\HH^{n}/\G_{0}$, and let $\T$ be the finite group coming from the previous lemma.

Now, we remark that the case of $n=2$ follows from the result of Long-MacLachlan-Reid \cite{LMR06}.
If $n=2$ then there is an index 2 subgroup $\G'<\G_{0}$
which is a genus 0 Fuchsian group. Thus, by Theorem 1.1 of \cite{LMR06}, we can conclude that
there are only finitely many possible $\G'$, since $\G'$ is also a congruence group, and therefore
there are only finitely many $\G_{0}$ and $\G$. In fact, $\Area(\HH^{2}/\G') < \frac{128}{3}\pi$.
Also, the case of $n=3$ was proven in \cite{Agol05}.
So we may restrict to the case that $4\leq n \leq 995$, by Prokhorov's Theorem \cite{Prokhorov86}.
In fact, since Nikulin has shown the finiteness of maximal non-uniform arithmetic reflection groups \cite{Nikulin80, Nikulin81, Nikulin87},
we may restrict to the case that $\G$ is cocompact, and therefore $n<30$, although our argument
gives a new proof of Nikulin's Theorem, making use of Prokhorov's Theorem. 

Consider $\HH^{n}\subset \SS^{n}$ embedded conformally as the upper hemisphere
of $\SS^{n}$, so that $\Isom(\HH^{n})$ acts conformally on $\SS^{n}$.
Normalize so that $\T$ acts isometrically on $\SS^{n}$, by placing the common fixed
point of $\T$ at the north pole of $\SS^{n}$. Clearly, $V_{PC}(n,\mathcal{O})=V_{PC}(n,P/\T)$,
since $|\OO|=|P/\T|$.
Then the orbifold $P/\T \subset \SS^{n}/\T$ is a conformal embedding, so fact 4 of section \ref{conformal volume} gives the inequality
$V_{PC}(n,P/\T) \leq V_{PC}(n,\SS^{n}/\T)$.

We now apply the eigenvalue estimates from Theorems \ref{ev} and \ref{BurgerSarnak}:
$$\frac{2n-3}{4} \Vol(\mathcal{O})^{\frac2n} \leq \l_{1}(\mathcal{O}) \Vol(\mathcal{O})^{\frac2n}  \leq n V_{PC}(n,\OO)^{\frac2n}\leq n (V_{PC}(\SS^{n}/\T))^{\frac2n}.$$

Assume now that $\G_{0}$ is defined over the number field $k$. By Theorem \ref{ellipticconformalvolume} and
Corollary  \ref{finitesubgroup}
we have 
$$V_{PC}(\SS^{n}/\T) \leq C(n) |\T|^{2n} \le C(n) (m_{n} \deg(k)^{2(n+1) n})^{2n}.$$
By  Lemmas \ref{volbydisc} and \ref{discriminant}, we
have $\Vol(\mathcal{O}) = \Vol(\HH^{n}/\G_{0}) \geq c_{n}c^{\delta_n \deg(k)}$. Putting these inequalities together, we obtain an inequality
of the form
$$c^{\delta_n\deg(k)} \leq C_{1}\!(n)\: \deg(k)^{4(n+1) n^{2}},$$
which implies that $\deg(k)$  must be bounded.

We also obtain
$$ \Vol(\mathcal{O})^{\frac2n} \le \frac{4n}{2n-3} C(n)^{\frac2n} m_{n}^4 \deg(k)^{8(n+1) n}.$$
Since $\deg(k)$ is bounded, this gives an upper bound on $\Vol(\OO)$ depending only on $n$.

Since volumes of arithmetic hyperbolic orbifolds are  discrete \cite{Borel81, Wang72}, and $\G$ is uniquely determined by $\G_{0}$,
we conclude that there are only finitely many arithmetic maximal reflection groups in dimension $n$.
\end{proof}

\section{Conclusion}
It is clear that for a finite volume polyhedron
$P\subset \HH^{n}$, $V_{PC}(n,P)= \Vol(\SS^{n})$, so by Theorem \ref{ev}
$$\l_{1}(P) \Vol(P)^{\frac2n} \le n \Vol(\SS^{n})^{\frac2n}.$$
Thus  in $n$ dimensions, there
are finitely many reflection groups $\G< \Isom(\HH^{n})$ which have a lower bound on
$\l_{1}(\HH^{n}/\G)$. If it were true that maximal arithmetic hyperbolic reflection groups
are congruence, then this would simplify the argument in this paper quite a bit, since
we would immediately get an upper bound on the volume by the above observation. So it
is natural to ask, does there exist a non-congruence maximal arithmetic reflection group?

The following conjecture seems to be interesting, and would also give an alternative argument
to the main theorem of this paper:
\begin{conjecture}
There is a function $K(n)$, such that if $\mathcal{O}$ is an elliptic $n$-orbifold,
then $V_{PC}(\mathcal{O})\leq K(n)$.
\end{conjecture}

Another approach to the proof of the main theorem would be possible if we assume the validity
of the {\it short geodesic conjecture} for arithmetic hyperbolic $n$-orbifolds.
As it was pointed out earlier in this section, Theorem \ref{ev} implies
$$\l_{1}(P) \Vol(P)^{\frac2n} \le n \Vol(\SS^{n})^{\frac2n}.$$
Now if the group $\Gamma$ generated by reflections in $P$ is arithmetic then it is contained in
some maximal arithmetic subgroup $\Gamma_0 < \Isom(\HH^{n})$ with a finite index $l$. Let $k$ be
the field of definition of $\Gamma_0$. The corollary from Lemma~\ref{cyclicvolume} and Lemma~\ref{vinb2}
imply that $l\le m_n\deg(k)^{2(n+1) n}$, and Lemmas~\ref{volbydisc}, \ref{discriminant} imply that
$\Vol(P) \ge lc^{\delta_n \deg(k)}$. What remains is to estimate $\l_{1}(P) = \l_{1}(\Gamma)$
using the fact that $\l_{1}(\Gamma_0) \ge (2n-3)/4$ which follows from the Theorem~\ref{BurgerSarnak}.
This brings us to the question about the spectrum of coverings which was studied by Brooks in
\cite{Brooks81, Brooks86}. The result of \cite{Brooks86} implies that
$\l_{1}(\Gamma) \ge c(\Gamma_0) \l_{1}(\Gamma_0)/l^2$ and more careful analysis
of the proof shows that the constant $c(\Gamma_0)$ depends only on the lower bound
for the injectivity radius of $\HH^n/\Gamma_0$. The short geodesic conjecture says
that in every dimension $n\ge 3$ there exist a universal lower bound for the injectivity
radius, so assuming the conjecture is true we can write $\l_{1}(\Gamma) \ge c(n) \l_{1}(\Gamma_0)/l^2$.
Substituting this to the main inequality together with the bounds for the index $l$
and $\Vol(P)$ we are getting the finiteness results in a similar way to the proof of
the theorem in sect.~6.
This argument avoids the analysis of the conformal volumes of spherical orbifolds
(sect.~3) but unfortunately is based on an open conjecture.

From the main theorem, we conclude that given an arithmetic reflection group  $\G<\Isom(\HH^{n})$, it must lie in
one of finitely many maximal reflection groups (up to conjugacy). If $\G$ is a reflection group in a polyhedron $P$
for which all the dihedral angles are $\pi/2$ or all are $\pi/3$, then
there are infinitely many finite covolume reflection subgroups of $\G$. Thus we see
that there are commensurability classes of arithmetic groups for which there are infinitely many
reflection groups in the commensurability class, and thus in our finiteness result, the
maximality assumption is crucial. The same is true for the arithmeticity assumption.
For example, consider the groups generated by reflections of the hyperbolic plane in the
sides of a hyperbolic triangle with angles $\pi/2$, $\pi/3$, $\pi/m$ ($m\ge 7$).
These groups are known to be maximal discrete subgroups of $\Isom(\HH^{2})$ but
all except finitely many of them are non-arithmetic. A similar construction is available
for  hyperbolic $3$-space with triangular prisms replacing triangles \cite{MR03}.
\medskip

It is an interesting project to try to identify all arithmetic maximal reflection groups
and to classify their reflection subgroups of finite covolume. Although all the constants
in our finiteness argument are computable, the quantitative estimates which could be deduced
directly from the proof are unreasonably large. However, if we restrict to the congruence
reflection groups (see the discussion in the beginning of this section) or make some additional
assumptions on the groups of automorphisms of the fundamental polyhedra of the reflection groups,
then the bounds can be improved essentially. 

\bibliographystyle{hamsplain}
\bibliography{higherarithmeticreflection5}

\def\cprime{$'$} \def\cprime{$'$} \def\cprime{$'$}
\providecommand{\bysame}{\leavevmode\hbox to3em{\hrulefill}\thinspace}
\providecommand{\href}[2]{#2}
\begin{thebibliography}{10}

\bibitem{Agol05}
I. Agol, \emph{{Finiteness of arithmetic Kleinian reflection groups}},
  Proceedings of ICM 2006, Madrid, vol.~II, European Mathematical Society,
  Zurich, 2006, \mbox{arXiv:math.GT/0512560}, pp.~951--960.

\bibitem{Allcock06}
D. Allcock, \emph{Infinitely many hyperbolic {C}oxeter groups through
  dimension 19}, Geom. Topol. \textbf{10} (2006), 737--758 (electronic).

\bibitem{Bel}
M.~Belolipetsky, \emph{Counting maximal arithmetic subgroups}, preprint,
  \mbox{arxiv:math.GR/0501198}.

\bibitem{Borel81}
A.~Borel, \emph{Commensurability classes and volumes of hyperbolic
  {$3$}-manifolds}, Ann. Scuola Norm. Sup. Pisa Cl. Sci. (4) \textbf{8} (1981),
  no.~1, 1--33.

\bibitem{BorelPrasad}
A. Borel and G. Prasad, \emph{Finiteness theorems for discrete subgroups
  of bounded covolume in semi-simple groups}, Inst. Hautes \'Etudes Sci. Publ.
  Math. (1989), no.~69, 119--171.

\bibitem{Brooks81}
R. Brooks, \emph{The fundamental group and the spectrum of the
  {L}aplacian}, Comment. Math. Helv. \textbf{56} (1981), no.~4, 581--598.

\bibitem{Brooks86}
\bysame, \emph{The spectral geometry of a tower of coverings}, J. Differential
  Geom. \textbf{23} (1986), no.~1, 97--107.

\bibitem{BurgerSarnak91}
M.~Burger and P.~Sarnak, \emph{Ramanujan duals. {II}}, Invent. Math.
  \textbf{106} (1991), no.~1, 1--11.

\bibitem{Dennin71}
J.~B. Dennin, Jr., \emph{Fields of modular functions of genus {$0$}},
  Illinois J. Math. \textbf{15} (1971), 442--455.

\bibitem{Dennin72}
\bysame, \emph{Subfields of {$K(2\sp{n})$} of genus {$0$}}, Illinois J. Math.
  \textbf{16} (1972), 502--518.

\bibitem{ElSoufiIlias86}
A.~El~Soufi and S.~Ilias, \emph{Immersions minimales, premi\`ere valeur propre
  du laplacien et volume conforme}, Math. Ann. \textbf{275} (1986), no.~2,
  257--267.

\bibitem{EGM98}
J.~Elstrodt, F.~Grunewald, and J.~Mennicke, \emph{Groups acting on hyperbolic
  space}, Springer Monographs in Mathematics, Springer-Verlag, Berlin, 1998,
  Harmonic analysis and number theory.

\bibitem{GJ78}
S. Gelbart and H. Jacquet, \emph{A relation between automorphic
  representations of {${\rm GL}(2)$} and {${\rm GL}(3)$}}, Ann. Sci. \'Ecole
  Norm. Sup. (4) \textbf{11} (1978), no.~4, 471--542.

\bibitem{GromovRohlin70}
M.~L. Gromov and V.~A. Rohlin, \emph{Imbeddings and immersions in {R}iemannian
  geometry}, Uspehi Mat. Nauk \textbf{25} (1970), no.~5 (155), 3--62.

\bibitem{Hersch70}
J. Hersch, \emph{Quatre propri\'et\'es isop\'erim\'etriques de membranes
  sph\'eriques homog\`enes}, C. R. Acad. Sci. Paris S\'er. A-B \textbf{270}
  (1970), A1645--A1648.

\bibitem{Horozov}
I.~Horozov, \emph{Euler characteristics of arithmetic groups}, Ph.D. thesis,
  Brown University, 2003, \mbox{arXiv:math.GR/0311117}.

\bibitem{Horozov05}
\bysame, \emph{Euler characteristics of arithmetic groups}, Math. Res. Lett.
  \textbf{12} (2005), no.~2-3, 275--291.

\bibitem{JL70}
H.~Jacquet and R.~P. Langlands, \emph{Automorphic forms on {${\rm GL}(2)$}},
  Springer-Verlag, Berlin, 1970, Lecture Notes in Mathematics, Vol. 114.

\bibitem{LangANT}
S. Lang, \emph{Algebraic number theory}, Addison-Wesley Publishing Co.,
  Inc., Reading, Mass.-London-Don Mills, Ont., 1970.

\bibitem{LiYau82}
P. Li and S.~T. Yau, \emph{A new conformal invariant and its
  applications to the {W}illmore conjecture and the first eigenvalue of compact
  surfaces}, Invent. Math. \textbf{69} (1982), no.~2, 269--291.

\bibitem{LMR06}
D.~D. Long, C.~Maclachlan, and A.~W. Reid, \emph{Arithmetic {F}uchsian groups
  of genus zero}, Pure Appl. Math. Q. \textbf{2} (2006), no.~2, 569--599.

\bibitem{MR03}
Colin Maclachlan and Alan~W. Reid, \emph{The arithmetic of hyperbolic
  3-manifolds}, Graduate Texts in Mathematics, vol. 219, Springer-Verlag, New
  York, 2003.

\bibitem{Margulis75}
G.~A. Margulis, \emph{Discrete groups of motions of manifolds of nonpositive
  curvature}, Proceedings of the International Congress of Mathematicians
  (Vancouver, B.C., 1974), Vol. 2, Canad. Math. Congress, Montreal, Que., 1975,
  pp.~21--34.

\bibitem{Nikulin80}
V.~V. Nikulin, \emph{On the arithmetic groups generated by reflections in
  {L}oba\v{c}evski\u{i} spaces}, Izv. Akad. Nauk SSSR Ser. Mat. \textbf{44}
  (1980), no.~3, 637--669, 719--720.

\bibitem{Nikulin81}
\bysame, \emph{On the classification of arithmetic groups generated by
  reflections in {L}oba\v{c}evski\u{i} spaces}, Izv. Akad. Nauk SSSR Ser. Mat.
  \textbf{45} (1981), no.~1, 113--142, 240.

\bibitem{Nikulin87}
\bysame, \emph{Discrete reflection groups in {L}obachevsky spaces and algebraic
  surfaces}, Proceedings of the International Congress of Mathematicians, Vol.
  1, 2 (Berkeley, Calif., 1986) (Providence, RI), Amer. Math. Soc., 1987,
  pp.~654--671.

\bibitem{Nikulin06}
V.~V. Nikulin, \emph{{Finiteness of the number of arithmetic groups generated
  by reflections in Lobachevsky spaces}}, preprint 2006,
  \mbox{arXiv:math.AG/0609256}, 5 pages.

\bibitem{Od90}
A.~M. Odlyzko, \emph{Bounds for discriminants and related estimates for class
  numbers, regulators and zeros of zeta functions: a survey of recent results},
  S\'em. Th\'eor. Nombres Bordeaux (2) \textbf{2} (1990), no.~1, 119--141.

\bibitem{PR}
V. Platonov and A. Rapinchuk, \emph{Algebraic groups and number
  theory}, Pure and Applied Mathematics, vol. 139, Academic Press Inc., Boston,
  MA, 1994, Translated from the 1991 Russian original by Rachel Rowen.

\bibitem{Prokhorov86}
M.~N. Prokhorov, \emph{Absence of discrete groups of reflections with a
  noncompact fundamental polyhedron of finite volume in a {L}oba\v{c}evski\u{i}
  space of high dimension}, Izv. Akad. Nauk SSSR Ser. Mat. \textbf{50} (1986),
  no.~2, 413--424.

\bibitem{ReedSimon78}
M. Reed and B. Simon, \emph{Methods of modern mathematical physics.
  {IV}. {A}nalysis of operators}, Academic Press [Harcourt Brace Jovanovich
  Publishers], New York, 1978.

\bibitem{ChernousovRyzhkov97}
A.~A. Ryzhkov and V.~I. Chernousov, \emph{On the classification of maximal
  arithmetic subgroups of simply connected groups}, Mat. Sb. \textbf{188}
  (1997), no.~9, 127--156.

\bibitem{S}
P. Scott, \emph{Subgroups of surface groups are almost geometric}, Journal
  of the London Mathematical Society. Second Series \textbf{17} (1978), no.~3,
  355--565.

\bibitem{Selberg65}
A. Selberg, \emph{On the estimation of {F}ourier coefficients of modular
  forms}, Proc. Sympos. Pure Math., Vol. VIII, Amer. Math. Soc., Providence,
  R.I., 1965, pp.~1--15.

\bibitem{Szego54}
G.~Szeg{\"o}, \emph{Inequalities for certain eigenvalues of a membrane of given
  area}, J. Rational Mech. Anal. \textbf{3} (1954), 343--356.

\bibitem{Thurston97}
W.~P. Thurston, \emph{Three-dimensional geometry and topology. {V}ol. 1},
  Princeton Mathematical Series, vol.~35, Princeton University Press,
  Princeton, NJ, 1997, Edited by Silvio Levy.

\bibitem{Tits66}
J.~Tits, \emph{Classification of algebraic semisimple groups}, Algebraic Groups
  and Discontinuous Subgroups (Proc. Sympos. Pure Math., Boulder, Colo., 1965),
  Amer. Math. Soc., Providence, R.I., 1966, 1966, pp.~33--62.

\bibitem{Vigneras83}
M.-F. Vign{\'e}ras, \emph{Quelques remarques sur la conjecture {$\lambda
  \sb{1}\geq {\frac14}$}}, Seminar on number theory, Paris 1981--82 (Paris,
  1981/1982), Progr. Math., vol.~38, Birkh\"auser Boston, Boston, MA, 1983,
  pp.~321--343.

\bibitem{Vinberg67}
{\`E}.~B. Vinberg, \emph{Discrete groups generated by reflections in
  {L}oba\v{c}evski\u{i} spaces}, Mat. Sb. (N.S.) \textbf{72 (114)} (1967),
  471--488; correction, ibid. 73 (115) (1967), 303.

\bibitem{Vinberg84}
\bysame, \emph{Absence of crystallographic groups of reflections in
  {L}oba\v{c}evski\u{i} spaces of large dimension}, Trudy Moskov. Mat. Obshch.
  \textbf{47} (1984), 68--102, 246.

\bibitem{Wang72}
H.~C. Wang, \emph{Topics on totally discontinuous groups}, Symmetric
  spaces (Short Courses, Washington Univ., St. Louis, Mo., 1969--1970), Dekker,
  New York, 1972, pp.~459--487. Pure and Appl. Math., Vol. 8.

\bibitem{YangYau80}
P.~C. Yang and S.~T. Yau, \emph{Eigenvalues of the {L}aplacian of
  compact {R}iemann surfaces and minimal submanifolds}, Ann. Scuola Norm. Sup.
  Pisa Cl. Sci. (4) \textbf{7} (1980), no.~1, 55--63.

\bibitem{Zograf91}
P.~Zograf, \emph{A spectral proof of {R}ademacher's conjecture for congruence
  subgroups of the modular group}, J. Reine Angew. Math. \textbf{414} (1991),
  113--116.

\end{thebibliography}

%    Bibliographies can be prepared with BibTeX using amsplain,
%    amsalpha, or (for "historical" overviews) natbib style.
%\bibliographystyle{amsplain}
%    Insert the bibliography data here.

\end{document}